\documentclass[pdflatex,sn-mathphys]{sn-jnl}% Math and Physical Sciences Reference Style
\usepackage{lineno,hyperref}
\usepackage{hyperref}       % hyperlinks
\usepackage{url}            % simple URL typesetting
\usepackage{booktabs}       % professional-quality tables
\usepackage{amsfonts}       % blackboard math symbols
\usepackage{nicefrac}       % compact symbols for 1/2, etc.
\usepackage{microtype}      % microtypography
\usepackage{siunitx}
\usepackage{enumitem}
\usepackage{algorithm}
\usepackage{bm}
\usepackage{amsmath,amssymb}
\usepackage{graphicx}
\jyear{2021}%

\theoremstyle{thmstyleone}%
\newtheorem{theorem}{Theorem}%  meant for continuous numbers
\newtheorem{lemma}[theorem]{Lemma}

\theoremstyle{thmstyletwo}%

\theoremstyle{thmstylethree}%

\begin{document}

\title[Functional Linear Partial Quantile Regression with Guaranteed Convergence]{Functional Linear Partial Quantile Regression with Guaranteed Convergence for Neuroimaging Data Analysis}

%%=============================================================%%
%% Prefix	-> \pfx{Dr}
%% GivenName	-> \fnm{Joergen W.}
%% Particle	-> \spfx{van der} -> surname prefix
%% FamilyName	-> \sur{Ploeg}
%% Suffix	-> \sfx{IV}
%% NatureName	-> \tanm{Poet Laureate} -> Title after name
%% Degrees	-> \dgr{MSc, PhD}
%% \author*[1,2]{\pfx{Dr} \fnm{Joergen W.} \spfx{van der} \sur{Ploeg} \sfx{IV} \tanm{Poet Laureate} 
%%                 \dgr{MSc, PhD}}\email{iauthor@gmail.com}
%%=============================================================%%

\author[1]{\fnm{Dengdeng} \sur{Yu}}\email{dengdeng.yu@uta.edu}

\author[2]{\fnm{Matthew} \sur{Pietrosanu}}\email{pietrosa@ualberta.ca}
\author[2]{\fnm{Ivan} \sur{Mizera}}\email{imizera@ualberta.ca}
\author[2]{\fnm{Bei} \sur{Jiang}}\email{bei1@ualberta.ca}
\author*[2]{\fnm{Linglong} \sur{Kong}}\email{lkong@ualberta.ca}
\author[3]{\fnm{Wei} \sur{Tu}}\email{wei.tu@queensu.ca}
%\equalcont{These authors contributed equally to this work.}

% \author[1,2]{\fnm{Third} \sur{Author}}\email{pietrosa@ualberta.ca}
% %\equalcont{These authors contributed equally to this work.}

\affil[1]{Department of Mathematics, University of Texas at Arlington, Arlington, TX 76019, USA}

\affil*[2]{Department of Mathematical and Statistical Sciences, 
  University of Alberta, 
   Edmonton, AB, T6G 2G1, Canada}

\affil[3]{Department of Public Health Sciences and Canadian Cancer Trials Group, Queen's University,  Kingston, ON, K7L 3N6, Canada}

%%==================================%%
%% sample for unstructured abstract %%
%%==================================%%

\abstract{Functional data such as curves and surfaces have become more and more common with modern technological advancements. The use of functional predictors remains challenging due to its inherent infinite-dimensionality. The common practice is to project functional data into a finite dimensional space. The popular partial least square (PLS) method has been well studied for the functional linear model \cite{delaigle2012methodology}. As an alternative, quantile regression provides a robust and more comprehensive picture of the conditional distribution of a response when it is non-normal, heavy-tailed, or contaminated by outliers. While partial quantile regression (PQR) was proposed in \cite{yu2016partial}, no theoretical guarantees were provided due to the iterative nature of the algorithm and the non-smoothness of quantile loss function. 
To address these issues,  we propose an alternative PQR (APQR) formulation with guaranteed convergence. This novel formulation motivates new theories and allows us to establish asymptotic properties. Numerical studies on a benchmark dataset show the superiority of our new approach. We also apply our novel method to a functional magnetic resonance imaging (fMRI) data to predict attention deficit hyperactivity disorder (ADHD) and a diffusion tensor imaging (DTI) dataset to predict Alzheimer's disease (AD). }

\keywords{Functional data analysis, Functional partial least square, Functional partial quantile regression, Finite smoothing, Neuroimaging}

\maketitle

\section{Introduction}

Functional data is becoming increasingly common in fields such as neuroimaging, finance, and other areas with longitudinal data, for example. While functional data carries rich information, incorporating them into statistical models brings extra challenges due to their inherent infinite-dimensionality. The common practice is to project functional data into a functional space with a finite functional basis. There are three major types of functional basis: general bases such as B-splines and wavelets, functional principal component (fPC) bases, and partial least square (PLS) bases. Relative to general bases, fPC bases incorporate variations in the predictors, and PLS bases consider variations in both the predictors and responses and are superior in supervised learning. The methodology and theory of PLS bases in functional linear models have been studied extensively in the literature \cite{delaigle2012methodology}. 
Linear models, however, often impose strong assumptions on the data-generating process. As an alternative to the linear model, quantile regression \cite{koenker1978regression} has been further developed and applied to functional linear regression \cite{kato2012estimation, yu2016partial}. It is well-known that, relative to least squares regression, quantile regression is more efficient and robust when the response is non-normal, when errors are heavy-tailed, or when outliers are present. Quantile regression is also capable of dealing with heteroscedasticity and provides a more complete picture of the dependence structure of the response \cite{koenker2005quantile}.

In this paper, we focus on the functional linear quantile regression model
\begin{eqnarray}
\label{chap2flqm}
Q_\tau(y \mid \mathbf{x}, \mathbf{z}(t)) = \alpha_\tau + \mathbf{x}^T \boldsymbol\beta_\tau + \int_I \mathbf{z}^T(t) \boldsymbol\gamma_\tau(t) dt,
\end{eqnarray}
where $Q_{\tau}(y \mid \mathbf x, \mathbf z(t))$ is the $\tau$-th conditional quantile of a response $y$ given scalar covariates $\mathbf x$ and functional covariates $ \mathbf z$, $\boldsymbol\beta_\tau$ is a vector of coefficients, and $\boldsymbol\gamma_\tau$ is a vector of functional coefficients, for a fixed quantile level $\tau \in (0,1)$. When estimating $\boldsymbol\gamma_\tau$, it is common to require that $\boldsymbol\gamma_\tau$ lies in a functional space and satisfies certain smoothness conditions \cite{james2009functional}.
Even with this restriction, however, $\boldsymbol\gamma_\tau$ generally still lies in an infinite-dimensional space and is difficult to estimate.

General bases \cite{cardot2005quantile,sun2005semiparametric} of B-splines \cite{cardot2003spline} or wavelets \cite{zhao2012wavelet}, as well as fPC bases, which use functional principal components to explain variation in $\mathbf{z}$ \cite{kato2012estimation,lu2014functional,tang2014partial,hall2007methodology,lee2012sparse}, are all well-documented in the existing literature \cite{cardot2005quantile,sun2005semiparametric}. In this article, we focus primarily on PLS bases and their estimation in the functional linear partial quantile regression model.

In traditional linear regression, PLS is an iterative procedure for feature selection originally developed for high-dimensional, multicollinear settings and is particularly popular in chemometrics \cite{wold1975soft,helland1990partial,frank1993statistical,nguyen2004partial,abdi2010partial}. PLS is a supervised dimension reduction technique that projects data into a lower-dimensional subspace formed by linear projections of covariates that are optimal for predicting a response. This method has been applied to functional data in \cite{preda2005pls} in the functional linear model, with consistency and convergence properties later established in \cite{delaigle2012methodology}. Analogous to the PLS basis in functional linear regression, the PLS basis in functional linear quantile regression is determined using information from both the response and functional covariates. This is not true of general and fPC bases which do not necessarily explain variation in the response. These two types of basis often require a large number of elements that may, in turn, contribute to model overfitting in the absence of regularization penalties \cite{cardot2005quantile,zhao2015wavelet}. On the other hand, a functional PLS basis can capture more relevant information with fewer terms \cite{delaigle2012methodology} and is thus an attractive alternative to general and fPC bases. Functional PLS bases in the functional linear quantile regression setting, which we refer to as partial quantile regression (PQR) bases, can be obtained using the PQR methodology and simple PQR (SIMPQR) algorithm proposed in \cite{yu2016partial} in which basis elements are extracted by sequentially maximizing the partial quantile covariance between the response and linear projections of functional covariates.

Despite the advantages of the PQR methodology, it is difficult to study its asymptotic properties due to its iterative nature and the non-differentiability of the quantile loss function. A reformulation of the functional PQR problem following the method of \cite{delaigle2012methodology} for traditional functional linear PLS is not possible in the quantile regression setting due to the non-additivity of conditional quantiles. Furthermore, most existing methods with ``nice'' properties, such as gradient descent, are no longer applicable due to the non-differentiability of the quantile loss function \cite{wu2009variable,zheng2011gradient}. Therefore an alternative method to PQR is necessary to achieve the guaranteed convergence of PQR. 
%A modified approach to PQR is therefore necessary in order to study the method's asymptotic properties.

% 
To address these problems, we propose an alternative PQR (APQR) formulation with guaranteed convergence. The contributions of the paper are three-folded. First, the proposed APQR formulation allows us to establish the guaranteed convergence of PLS basis in functional quantile regression model, and various asymptotic properties such as consistency and asymptotic normality for APQR estimators are established. Second. the proposed formulation has two main components, a smooth approximation of the quantile loss function $\rho_\tau$ using finite smoothing techniques \cite{muggeo2012quantile,chen2012finite,crambes2009smoothing,zhao2015wavelet, 8970747, zhu2017expectile} and a simultaneous estimation of the PQR basis. %using a block relaxation method \cite{de1994block}. 
Third. applications to benchmark dataset and large scale neuroimaging datasets show the superiority and usefulness of the proposed method.

The rest of this paper is organized as follows. 
In Section \ref{chap2sec3}, we propose the APQR method, 
and define the APQR for extracting PQR bases in the functional linear quantile regression model.
A number of asymptotic properties of APQR are discussed in Section \ref{chap2sec4}, including score, information derivations, and model identifiability properties, estimation consistency, asymptotic normality, and guaranteed convergence. 
Detailed proofs and verifications for the key assumptions are provided in the appendix following the main text.
The algorithm is implemented in Section \ref{APQR algorithm}.  
In Section \ref{chap2sec5}, APQR is applied to simulated and real biomedical datasets collected by the Alzheimer's Disease Neuroimaging Initiative (ADNI) and ADHD-200 from the International Neuroimaging Datasharing Initiative (INDI) to demonstrate its superiority or comparability over SIMPQR and other existing methods.

\section{Alternative Partial Quantile Regression}
\label{chap2sec3}

In our new APQR method, we propose estimating model (\ref{chap2flqm}), for given basis size $K$, by solving
\begin{eqnarray*}
\label{chap2flqm2}
\min_{\alpha,\boldsymbol\beta,\phi_{\tau 1},\dots,\phi_{\tau K}}
\mathbb{E}\ \rho_\tau \left( Y - \alpha - \mathbf{X}^T\boldsymbol\beta - \sum_{k=1}^K\int_0^1 Z(t) \phi_{\tau k}(t) dt \right),
\end{eqnarray*}
that is, we obtain basis elements simultaneously rather than sequentially. Similar to the PLS formulation in a traditional functional linear model, we find a group of directions $\phi_{\tau 1},\dots,\phi_{\tau K}$ so that the projections $Z_1,\dots,Z_K$ of $Z$ onto these directions contribute as much as possible to predicting the level-$\tau$ conditional quantile of the response, after accounting for the other covariates $\mathbf{X}$. In other words, we obtain a PQR basis $\phi_{\tau 1},\dots,\phi_{\tau K}$ by maximizing $\mathrm{Cov}_{\tau} \left( Y,~ \sum_{k=1}^K \int_0^1 {Z}(t) \phi_{\tau k}(t) dt\ \Big\vert\ \mathbf{X} \right)$.

Despite the advantages of PQR basis, it is difficult to study the asymptotic properties of SIMPQR estimates due to the iterative nature of our previous approach \citep{yu2016partial} and the non-differentiability of the quantile loss function. This difficulty is also present in traditional PLS for functional linear regression. \cite{delaigle2012methodology} addressed this issue in the latter case by proposing an alternative but equivalent PLS formulation, called APLS, that exploited an equivalence between the fPC and functional APLS spaces in order to establish APLS consistency and convergence rates. Unfortunately, a similar equivalence does not exist between PQR and fPC spaces due to the non-additivity of conditional quantiles. In the following subsections, we explore modifications to address these issues, resulting in our new APQR approach and algorithm for PQR basis estimation. This approach motivates new theory and yields insightful results regarding asymptotic properties of this method through advanced techniques in empirical processes theory \cite{van2000asymptotic}.

Consider the functional linear quantile regression model (\ref{chap2flqm}). Fix some $\tau \in (0,1)$
and let $K$ be the dimension of the functional space into which $\boldsymbol\gamma_\tau$ is projected. We use $\left\{ \phi_{\tau k}(t)\ \lvert\ k=1,\dots,K\right\}$ to denote the PQR basis for this space. For each $i = 1,\ldots,n$, the response $Y$ is observed as $y_i$, while $Z(t)$ is observed only at $d$ discrete points in $\mathcal{T} = \{0=t_1 < t_2 < \ldots < t_d=1\}$. In other words, we only observe the functional response as $z_{i}(t_j)$ for $i=1,\dots,n$ and $j=1,\ldots,d$.
Define $\mathbf{Z} = (\mathbf{z}_{1},\ldots,\mathbf{z}_{n})^T
\in \mathbb{R}^{n \times d}$, where $\mathbf{z}_i = (z_i(t_1),\ldots,z_i(t_d))^T \in \mathbb{R}^d$ for $i=1,\ldots,n$. The PQR basis $\{\phi_{\tau k}\}_{k=1}^K$ can be obtained via $\mathbf{C}_\tau \in \mathbb{R}^{d \times K}$ by maximizing
\begin{eqnarray}
\label{chap2scova2}
l(\alpha,\boldsymbol\beta,\mathbf{C}) = -  \sum_{i=1}^{n} \rho_{\tau}\left(y_i - \alpha - \mathbf{x}_i^T \boldsymbol\beta-
\mathbf{z}_{i}^T \mathbf{C}_{} \mathbf{1}_K \right),
\end{eqnarray}
where $\mathbf{C}_\tau=(\mathbf{c}_{\tau 1},\ldots,\mathbf{c}_{\tau K}) \in \mathbb{R}^{d \times K}$ with
$\mathbf{c}_{\tau k} = (\phi_{\tau k}(t))_{t \in \mathcal{T}}\in \mathbb{R}^{d}$ for $k=1,\ldots,K$. Each column of $\mathbf{C}_\tau$ is a vector of functional basis values evaluated on the discrete domain $\mathcal{T}$. For a sequence of smoothing parameters $\left(\nu_N\right)_{N=1}^{\infty}$ converging to $\nu_0$, we use a uniformly smooth approximating function $\rho_{\tau \nu_N}$ that converges uniformly to $\rho_\tau$ as $N \to \infty$.

\subsection{Convex smoothing approximation}
The first hurdle in studying the asymptotic properties of SIMPQR basis estimates appears in the non-differentiability of the quantile loss function $\rho_\tau$ at zero. We propose a convex, finite smoothing approximation $\rho_{\tau \boldsymbol\nu}$ of $\rho_\tau$, where $\boldsymbol\nu$ is a vector of smoothing parameters, such that $\rho_{\tau \boldsymbol\nu}$ is differentiable and converges uniformly to $\rho_{\tau}$ as $\boldsymbol\nu \to \boldsymbol\nu_0$ \cite{muggeo2012quantile,chen2012finite}.
Because this approximation is convex and converges uniformly, we have that the minimizer of $\rho_{\tau \boldsymbol\nu}$ converges to the minimizer of $\rho_{\tau}$ on a compact set \cite{hjort2011asymptotics}.
After replacing $\rho_\tau$ with $\rho_{\tau\boldsymbol\nu}$, the quantile loss function becomes differentiable, making it possible to define score, information, and identifiability for functional linear quantile regression. This approach is an example of finite smoothing techniques applied to non-smooth convex optimization problems and is crucial to quantile regression, both statistically and computationally \cite{nesterov2005smooth,wu2009variable,zheng2011gradient,chen2012finite}. Various options are available for the smoothing function, including the \emph{generalized Huber function} \cite{chen2012finite}
and \emph{iterative least squares smoothing function} \cite{muggeo2012quantile}. For a given precision $\varepsilon$ and an appropriate smoothing function, one may obtain an optimal efficiency estimate in terms of the number of iterations until convergence as $O(\varepsilon^{-1})$.
This is a significant improvement over other popular numerical schemes for non-smooth convex minimization such as subgradient methods,
whose number of iterations until convergence is $O(\varepsilon^{-2})$
\cite{nesterov2005smooth}.

Replacing $\rho_\tau$ by $\rho_{\tau \nu_N}$ in (\ref{chap2scova2}), we approximate $l$ with
\begin{eqnarray}
\label{chap2scova3}
l_N(\alpha,\boldsymbol\beta,\mathbf{C}) = - \sum_{i=1}^{n} \rho_{\tau \nu_N}\left(y_i - \alpha - \mathbf{x}_i^T \boldsymbol\beta-
\mathbf{z}_{i}^T \mathbf{C}_{} \mathbf{1}_K \right),
\end{eqnarray}
where $l_N$ converges uniformly to $l$ as $N \to \infty$. A crucial observation here is that both $l$ and $l_N$ are blockwise convex about $(\alpha,\boldsymbol\beta)$ and $\mathbf C$. Our APQR algorithm for prespecified $K$ and $\tau$ follows below.

%\subsection{Empirical Choice of Basis Size}

%\textbf{Basis size selection} 
\subsection{Basis size selection}
As discussed previously, it is important to choose a PQR basis size $K$ large enough to adequately approximate Model (\ref{chap2flqm2}), but not so large to induce model overfitting. Numerous criteria are available for choosing $K$ such as CV prediction error or BIC, similar to the criteria adopted by other authors to choose an fPC basis size \cite{kato2012estimation,lu2014functional,tang2014partial}. In our case, we define CV error as $
\rm{CV}_\tau(K) =  \frac{1}{n}\sum_{i=1}^n \rho_\tau \left( y_i -\hat{\alpha}^{(-i)}_\tau - \mathbf{x}_i^T \hat{\boldsymbol\beta}^{(-i)}_\tau -
\sum_{k=1}^{K} \hat{\mathbf{z}}_{ki} \hat{\boldsymbol\gamma}^{(-i)}_{\tau k} \right)$, where $\hat{\boldsymbol\gamma}^{(-i)}_{\tau k}$ and $\hat{\alpha}^{(-i)}$ and $ \hat{\boldsymbol\beta}^{(-i)}_\tau$, for $i=1,\dots,n$ and $k=1,\ldots,K$, are estimates of Model (\ref{chap2flqm2}) parameters computed after removing the $i$-th observation.

\section{Asymptotic Properties}
\label{chap2sec4}
In the following section, we adopt the asymptotic setup with fixed $K$ and divergent sample size $n$, and further simplify notation by suppressing the subscript $\tau$. Since our focus is on the basis estimate $\hat{\mathbf{C}}$, we omit $\alpha$ and $\beta$, although the conclusions here generalize easily to include both. For a given, fixed $K$, we establish score, information, and identifiability. Consistency and asymptotic normality can be derived by applying empirical process theory \cite{van2000asymptotic} and by following discussions similar to those in \cite{li2013tucker} and \cite{zhou2013tensor}.

\begin{theorem}
\label{chap2prop1}
$l_{N}(\boldsymbol\upsilon)$ converges uniformly to $l(\boldsymbol\upsilon)$ with $\boldsymbol\upsilon=(\alpha,\boldsymbol\beta,\mathbf{C})$ as $N \to \infty$,
where $l_N(\boldsymbol\upsilon)$ and $l(\boldsymbol\upsilon)$ are as defined in (\ref{chap2scova2}) and (\ref{chap2scova3}).
For a fixed $N$, if (i) $l_{N}(\boldsymbol\upsilon)$ is continuous and coercive, that is, if the set $\{ \boldsymbol\upsilon : l_{N}(\boldsymbol\upsilon)\geq l_{N}(\boldsymbol\upsilon^{(0)})\}$ is compact and bounded above, and 
(ii) the objective function in each block update of the algorithm is strictly concave, and (iii) the set of stationary points of $l_{N}(\boldsymbol\upsilon)$ are isolated, we have that
\begin{enumerate}
\item (Global convergence) The sequence $\boldsymbol\upsilon^{(p)} = (\alpha^{(p)},\boldsymbol\beta^{(p)},\mathbf{C}^{(p)})$ generated by the APQR algorithm converges to a stationary point of $l_{N}(\boldsymbol\upsilon)$.
\item (Local convergence) Let $\boldsymbol\upsilon^{(\infty)} = (\alpha^{(\infty)},\boldsymbol\beta^{(\infty)},\mathbf{C}^{(\infty)})$ be a strict local maximum of $l_{N}(\boldsymbol\upsilon)$.
The iterates generated by the APQR algorithm above are locally attracted to $\boldsymbol\upsilon^{(\infty)}$ for $\boldsymbol\upsilon^{(0)}$ sufficiently close to $\boldsymbol\upsilon^{(\infty)}$.
\item (Approximation convergence) The convergence points obtained from $l_{N}(\boldsymbol\upsilon)$ converge in probability to the convergence point of $l(\boldsymbol\upsilon)$ as $N \to \infty$.
\end{enumerate}
\end{theorem}

The assumptions above are not hard to verify if we impose some regularity conditions on the distribution functions (\cite{koenker2005quantile}). If assumptions (i), (ii), and (iii) hold, then local and global convergence can be obtained following the discussions in \cite{li2013tucker} and \cite{zhou2013tensor}. In order to prove approximation convergence, provided that $l_N(\boldsymbol\upsilon)$ has a unique maximizer, we can apply Lemma 2 of \cite{hjort2011asymptotics}. The details verifying this result are deferred to the appendix.

%\subsection{Score and Information}
We then derive the score and information matrix for the functional linear quantile regression model in APQR. As discussed in \cite{yu2001bayesian} and \cite{sanchez2013likelihood}, minimizing the quantile loss is equivalent to maximizing the likelihood function formed by independent and identically distributed asymmetric Laplace densities. In fact, the function $l$ defined in Equation (\ref{chap2scova2}) is proportional to this log-likelihood. We seek to derive the score and information matrix using $l$ instead. However, since $l$ is not differentiable, we use its smoothing approximation $l_N$, defined in Equation (\ref{chap2scova3}). The difference in score and information matrices between $l_N$ and $l$ becomes negligible as $N \to \infty$.

\begin{theorem}[Score and information]
\label{chap2prop2}
Consider $l_N$ in Equation (\ref{chap2scova3}), as derived from the functional linear quantile regression model. Let $\eta_i(\mathbf{C}) = \mathbf{z}_{i}^T \mathbf{C} \mathbf{1}_K$ with gradient $\nabla \eta_i(\mathbf{C}) = - \mathbf{1}_K \otimes \rm{vec}\left(\mathbf{z}_{i} \right)\in \mathbb{R}^{K d}$.
\begin{enumerate}
\item The score function (or score vector) of $l_N$ is $\nabla l_N(\mathbf{C}) =  - \sum_{i=1}^{n} \rho^\prime_{\tau \nu_N}\left( \eta_i(\mathbf{C}) \right) \cdot \nabla \eta_i(\mathbf{C})$.
\item The Fisher information matrix of $l_N$ is given by \\
$\mathbf{I}_N (\mathbf{C}) =\sum_{i=1}^n \sum_{j=1}^n \rho^\prime_{\tau \nu_N} \left( \eta_i(\mathbf{C}) \right)
\rho^\prime_{\tau \nu_N} \left( \eta_j(\mathbf{C}) \right) $ $\nabla \eta_i(\mathbf{C}) d  \eta_j(\mathbf{C})$.
\end{enumerate}
\end{theorem}

%\subsection{Identifiability}
Before continuing to asymptotic properties, we address the issue of model identifiability in APQR. Parametrization of our functional linear quantile regression model is currently non-identifiable due to the indeterminacy of $\mathbf{C}$ up to column permutation. In other words, $\mathbf{C}$ and $\mathbf{C} \boldsymbol\Pi$ yield the same model for any $K \times K$ permutation matrix $\boldsymbol\Pi$. To address this indeterminacy, we assume that the elements comprising the first row of $\mathbf{C}=(c_{ij})_{ij}$ are distinct and arranged in decreasing order. Observe that the resulting parameter space $S=\{\mathbf{C}=(c_{ij})_{ij}\ \lvert\  c_{11} > \cdots > c_{1K}\}\subset \mathbb{R}^{d \times K}$ is open and convex. Having addressed identifiability, we next give a necessary and sufficient condition for local identifiability. The following result can be easily verified with Theorem 1 of \cite{rothenberg1971identification}. The details verifying this result are deferred to the appendix.

\begin{theorem}[Local identifiability]
Let $\{ (y_i,\mathbf{z}_{i})\}_{i=1}^n$ be a sequence of i.i.d. data points from the functional linear quantile regression model and let $\mathbf{C}_0 \in S$ be a point with a neighbourhood in which the information matrix $\mathbf{I}_N (\mathbf{C})$ has constant rank. Then $\mathbf{C}_0$ is identifiable up to permutation if and only if
$I_N(\mathbf{C}_0) = \sum_{i=1}^n \sum_{j=1}^n \rho^\prime_{\tau \nu_N} \left( \eta_i(\mathbf{C}_0) \right)
\rho^\prime_{\tau \nu_N} \left( \eta_j(\mathbf{C}_0) \right) \nabla \eta_i (\mathbf{C}_0) d  \eta_j(\mathbf{C}_0)
$
is nonsingular.
\end{theorem}

Asymptotic properties of our estimators follow from discussions on maximum likelihood estimators and M-estimation. A key observation here is that, by adopting the standard tensor notation from \cite{zhou2013tensor}, Model (\ref{chap2flqm}) can be rewritten as $Y = \alpha + \boldsymbol\beta^T \mathbf{X} + \langle \mathbf{C}\mathbf{1}_K, \mathbf{Z}\rangle$. It is known that the collection of degree 1 polynomials $\{\langle \mathbf{C} \mathbf{1}_K,\mathbf{Z}\rangle\ \lvert\ \mathbf{C} \in S \}$ form a Vapnik-Cervonenkis class. As a result, the standard uniform convergence theory for M-estimation applies \cite{van2000asymptotic}. 

\begin{theorem}[Consistency]
\label{chap2thmconsistency}
Assume that $\mathbf{C}_0 \in S \subset \mathbb{R}^{ d \times K}$
 is (globally) identifiable up to permutation and that the functional covariate vectors $\mathbf{Z}_i$ are i.i.d. from a bounded distribution. The M-estimator is consistent, that is, $\hat{\mathbf{C}}_n$ converges to $\mathbf{C}_0$ in probability (up to permutation) for the functional linear quantile regression model on a compact set $S_0 \subset S$.
\end{theorem}

Consistency can be verified using the theory of empirical processes. Further details are available in the appendix following this main text. By showing that $\left\{ l_N(\mathbf{C})\ \vert\ \mathbf{C} \in S\right\}$ is a Donsker class when parameters are restricted to a compact set, the Glivenko-Cantelli Theorem establishes uniform convergence. Uniqueness is guaranteed by the information equality whenever $\mathbf{C}_0$ is identifiable.

\begin{theorem}[Asymptotic normality]
\label{chap2thmnormality}
For an interior point
$\mathbf{C}_0 \in S$ with nonsingular information matrix $I_N \left(\mathbf{C}_n \right)$ and consistent estimator $\hat{\mathbf{C}}_n$, we have that $\sqrt{n} \left[\rm{vec}\left(\hat{\mathbf{C}}_n \right) - \rm{vec} \left(\mathbf{C}_0\right) \right]$ converges in distribution to a normal random variable with mean zero and covariance $I^{-1}_N \left(\mathbf{C}_0 \right)$.
\end{theorem}

\section{APQR algorithm}
\label{APQR algorithm}
We next discuss our modifications to the original SIMPQR algorithm proposed by \cite{yu2016partial}. In particular, we incorporate a block relaxation method \cite{de1994block} and take advantage of a back-fitting scheme to obtain and update PQR basis simultaneously instead of sequentially. This modification provides an alternative PQR (APQR) formation with asymptotic properties that can be obtained more easily than in SIMPQR. Throughout this subsection, we assume a fixed $K$.

\begin{algorithm*}[h]
\caption{Alternative SIMPQR (APQR), given $K$ and $\tau$}
\begin{enumerate}
\item {\it Initialize}:
Standardize $z_i(t_j)$ for each $j$ to have mean zero and variance one.
\item{\it Repeat} (for iteration number $N=1,2,\dots$):
\begin{enumerate}[labelindent=0.5cm]
\item {\it Initialize}: Let $\mathbf{C}^{(0)}_\tau$ be a random matrix (or otherwise initially specified),
with $(\alpha^{(0)}_\tau, \boldsymbol\beta^{(0)}_\tau) = \arg\max_{\alpha, \boldsymbol\beta} l_{N}(\alpha,\boldsymbol\beta,\mathbf{C}^{(0)}),$
\item {\it Repeat} (for update number $p=1,2,\dots$):
\begin{enumerate}
\item {\it Update}: 
$\mathbf{C}_\tau^{(p+1)} = \arg\max_{\mathbf{C}} l_{N}(\alpha_\tau^{(p)},\boldsymbol\beta_\tau^{(p)},\mathbf{C})$,
\item {\it Update}: 
$(\alpha^{(p+1)}_\tau, \boldsymbol\beta^{(p+1)}_\tau) $  $=\arg\max_{\alpha, \boldsymbol\beta} l_{N}(\alpha,\boldsymbol\beta,\mathbf{C}_\tau^{(p+1)}).$
\item {\it Save}:  $\alpha_\tau = \alpha_\tau^{(p+1)}$, $\boldsymbol\beta_\tau=\boldsymbol\beta^{(p+1)}_\tau$  and $\mathbf{C}_\tau = \mathbf{C}_\tau^{(p+1)} $.
\item {\it Stop} (for $p$) when \\
$\left\lvert l_{N} (\alpha^{(p+1)}_\tau, \boldsymbol\beta^{(p+1)}_\tau, \mathbf{C}_\tau^{(p+1)}) \right. $ $\left.- l_{N} (\alpha^{(p)}_\tau, \boldsymbol\beta^{(p)}_\tau, \mathbf{C}_\tau^{(p)}) \right\rvert < \varepsilon$.
\end{enumerate}
\end{enumerate}
\item {\it Stop} (for $N$) when 
$l_{N+1} (\alpha_\tau, \boldsymbol\beta_\tau,\mathbf{C}_\tau) - l_{N} (\alpha_\tau, \boldsymbol\beta_\tau,\mathbf{C}_\tau) < \epsilon$.
\item {\it Output}: $\alpha_\tau$, $\boldsymbol\beta_\tau$ and $\mathbf{C}_\tau $.
\item {\it Define}: $\tilde{\mathbf{Z}} =\mathbf{Z}\mathbf{C}_\tau = (\tilde{\mathbf{z}}_1,\ldots,\tilde{\mathbf{z}}_K)\in \mathbb{R}^{n \times K}$ as the projections of the $z_i$ onto the basis $\{\phi_{\tau k}\}_{k=1}^K$.
\item {\it Define}: $\hat{Q}_\tau(Y\mid\tilde{\mathbf{Z}})$ as the estimated conditional quantile.
\end{enumerate}
\end{algorithm*}

The advantages of APQR over SIMPQR are two-fold. First, the columns of $\mathbf{C}_\tau$ in APQR, which contain the PQR basis vectors, are all updated simultaneously using the back-fitting scheme, unlike the sequential updating scheme of SIMPQR.
Second, this block relaxation procedure guarantees solution stability and convergence because the objective function $l_N$ is blockwise convex and is bounded above. This fact is easy to verify.
Last but not least, the convexity and uniform convergence of the objective functions guarantee that
the maximizer of $l_N$ converges to the maximizer of $l$, thus assuring that the $\mathbf C_\tau$ obtained when optimizing $l_N$ is close to the true $\mathbf C_\tau$. Details have been given in Theorem \ref{chap2prop1} of the previous section.

\section{Numerical Studies} \label{chap2sec5}

We now study the finite sample performance of our proposed APQR method using one simulated and two real datasets. In all of our analyses, APQR shows superior or at least comparable performance relative to the three other methods considered.

\subsection{Benchmark phoneme dataset}
Our first analysis follows the setup of \cite{delaigle2012methodology}. We obtain functional covariate observations from the benchmark Phoneme dataset (available online at \url{web.stanford.edu/~hastie/ElemStatLearn/datasets/}), of which we only consider the data associated with "aa" and "ao" phonemes, yielding $1717$ instances, each with log-periodogram measurements for $d=256$ uniformly-spaced frequencies. The response for this simulation is generated according to the linear functional model $Y= \int_0^1 Z(t) \gamma(t) dt + \varepsilon$.
Here, the error $\varepsilon$ is taken to be Gaussian with mean zero and variance five times that of the observed data.  
We consider four different models through different specifications of the functional coefficient $\gamma$ by computing the first $J=20$ sample fPC basis functions $\hat\phi_j$ for $j=1,\dots,J$, and defining four different $\gamma(t) = \sum_{j=1}^Ja_j\hat\phi_j(t)$ by taking $\text{(i)}~a_j = (-1)^j\mathbf{1}_{1\leq j \leq 5},\  \text{(ii)}~a_j = (-1)^j\mathbf{1}_{6\leq j \leq 10},\  \text{(iii)}~a_j = (-1)^j\mathbf{1}_{11\leq j \leq 15},\  \text{(iv)}~a_j = (-1)^j\mathbf{1}_{16\leq j \leq 20}.$
It is clear that the above models are in order of decreasing favorability towards quantile regression with an fPC basis (QRfPC). On the other hand, we will see that the PLS, PQR, and APQR methods are able to extract the relationship between the functional covariates and the response with a relatively small number of basis elements.

\begin{figure*}
  \centering
 \includegraphics[width=\textwidth]{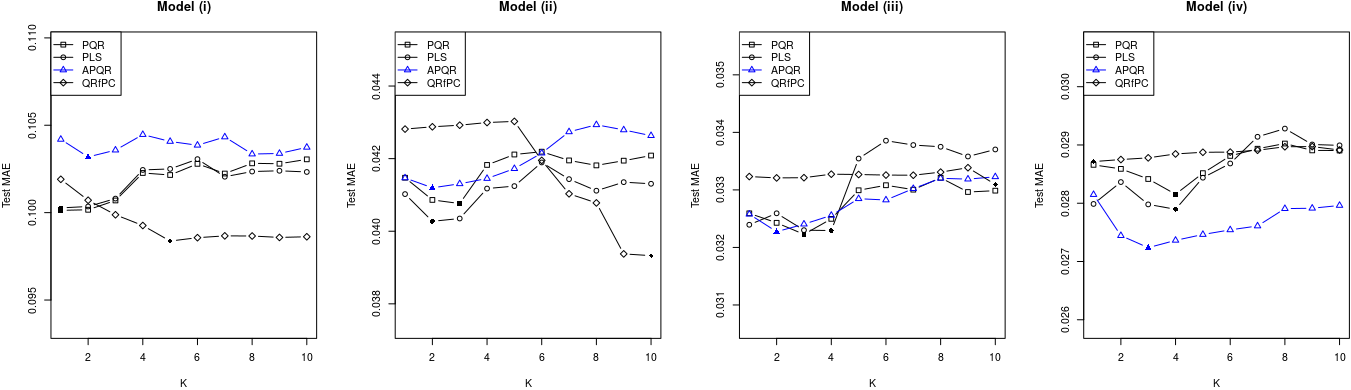}

 \caption{Test mean absolute error (MAE) averaged over 10 datasets vs. basis size $K$ for the four models. Each line represents one of the PQR, PLS, APQR, and QRfPC methods, with our proposed APQR method further highlighted in blue. A solid symbol denotes optimal basis size for each method, corresponding to minimal 10-fold cross-validation MAE.}

 \label{chap2sim2f1}
\end{figure*}

We generate $10$ datasets according to the above models and randomly split each dataset into training and tests sets composed of $1500$ and $217$ observations. We compute mean absolute (prediction) error (MAE) for the test set, averaged across all 10 simulated datasets. As $\tau=0.5$, MAE is equivalent to mean quantile error. We determine an optimal basis size $K$ for each method by minimizing 10-fold cross-validation MAE on the training set, averaged over all 10 datasets. A visual summary of these results is given in Figure \ref{chap2sim2f1}. As expected, QRfPC performs optimally in Model (i) where the functional coefficient $\gamma$ is in fact a linear combination of the first five sample fPC basis functions. This also explains the sudden decrease in test MAE for QRfPC in Model (ii) when $K$ increases from 5 to 10. APQR underperforms in Model (i) relative to the other three methods, but has increasingly better performance relative to other methods in Models (ii) through (iv), ultimately outperforming PQR, PLS, and QRfPC for $K=1,\dots,10$. In all models, PQR and PLS perform comparably.

\subsection{ADHD-200 fMRI dataset} 
We apply our proposed APQR method to an attention deficit hyperactivity disorder (ADHD) dataset available from New York University's Child Study Center as part of the ADHD-200 Sample Initiative Project (\url{http://fcon_1000.projects.nitrc.org/indi/adhd200/}). The response of interest is the revised Conners' parent rating scale (long version), a continuous ADHD severity index. Relevant scalar covariates include sex, age, handedness, ADHD diagnosis status, medication status, as well as verbal, performance, and Full4 IQ test scores. We separately consider six functional covariates, specifically, average grayscale values obtained from resting-state fMRI measurements of the cerebellum, vermis, and temporal, parietal, occipital, and frontal lobes. Averages are taken over measurements from a total of 116 regions of interest, each observed at $d=172$ equally-spaced time points. After data cleaning, we include $n=120$ subjects in our analysis with no missing value.

\begin{figure*}
 \centering
 \vspace{.3in}
%   \vspace{-1.0cm}
 \includegraphics[width = \textwidth]{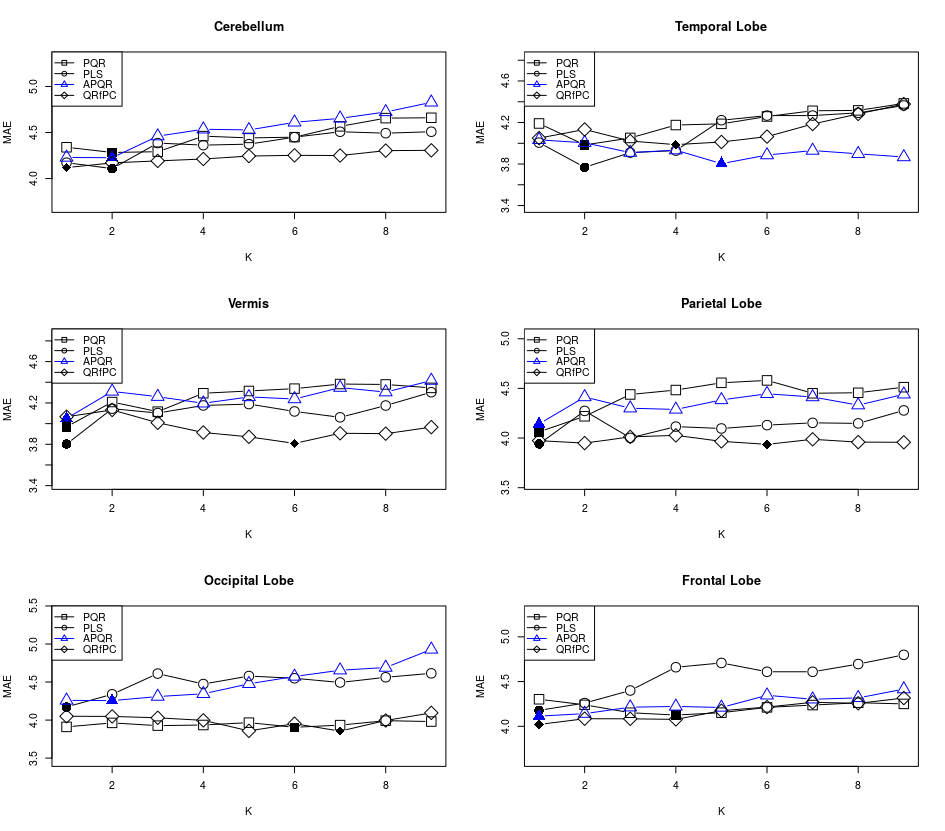}
 \vspace{.3in}
 % \vspace{-0.5cm}
 \caption{ Test mean absolute error (MAE) (averaged over 10 simulations) vs. basis size $K$ for the ADHD-200 fMRI data for six different functional covariates. Each line represents one of the PQR, PLS, APQR, and QRfPC methods, with our proposed APRQ method highlighted in blue. A solid symbol denotes optimal basis size for each method, as determined by 10-fold cross-validation MAE.}
 \label{chap2sim4f1}
\end{figure*}

We randomly split the dataset into training and test sets composed of 100 and 20 subjects, respectively, in 10 different ways. Optimal basis sizes $K$ are selected with average performance compared for each method as in the previous subsection. As illustrated in Figure \ref{chap2sim4f1}, APQR performs comparably to other methods in terms of both MAE for a given basis size $K$ and MAE at an optimal $K$. Excepting the case using functional data from the temporal lobe, we see that APQR has an optimal basis size of at most two. This is an improvement over the PQR and QRfPC methods, which generally requires at least four basis elements.

\begin{figure*}
\centering
% \vspace{-0.5cm}
\includegraphics[width = \textwidth]{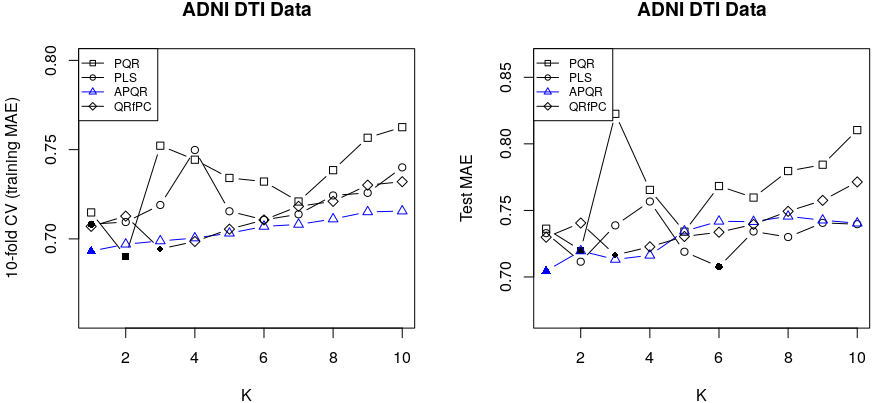}
 %\vspace{-0.5cm}
\caption{ 10-fold cross-validation (left) and test (right) mean absolute error (MAE) averaged over 10 datasets vs. basis size $K$ for the ADNI dataset. Each line represents one of the PQR, PLS, APQR, and QRfPC methods, with our proposed APQR method further highlighted in blue. Shaded symbols denote the optimal basis size for each method, as determined by 10-fold cross-validation MAE.}

\label{chap2sim5f1}
\end{figure*}
%\textbf{ADNI data}
\subsection{ADNI data}
We next consider a diffusion tensor imaging (DTI) dataset collected from the NIH Alzheimer's Disease Neuroimaging Initiative (ADNI) study (\url{adni.loni.ucla.edu}). The goal is to predict scores on the Mini Mental State Examination (MMSE) --- a widely-used Alzheimer's screening test. Scalar covariates for the ADNI dataset include sex, age, years of education, Alzheimer's diagnosis status, and APOE-4 allele genotype. We take fractional anisotropy (FA) values along the corpus callosum fiber tract as the single functional covariate. For each of the $n=200$ subjects included in our analysis, FA is measured at $d=83$ equally-spaced grid points and may be considered a function of arc length. The data have been previously analyzed in \cite{yu2019sparse, wang2019wavelet}.

We randomly split the dataset into training and test sets composed of 180 and 20 subjects, respectively, in 10 different ways. 
Figure \ref{chap2sim5f1} suggests that APQR and PLS generally outperform PQR and QRfPC, so that APQR is the best among the three quantile regression methods we considered for the ADNI dataset. For any given basis size $K$, APQR and PLS perform comparably to, if not better than, the other two methods. More notably, we see that APQR has an optimal basis size of $K=1$, while PLS, PQR, and QRfPC require $K=6$, $K=2$, and $K=3$, respectively. This result, consistent with our previous analyses, suggests that APQR better captures relationships between the functional covariate and the response with fewer basis elements.

\section{Conclusion}

Functional data such as curves and surfaces have become increasingly common with modern technological advancements. This type of data is often treated as discrete data points: generally, extracting useful knowledge from such data has been challenging due to its infinite dimensionality. 
While the popular partial least square (PLS) method has been well studied in the context of the linear functional model \cite{delaigle2012methodology}, quantile regression alternatively provides a robust and more comprehensive picture of the conditional distribution of the response when it is non-normal, heavy-tailed, or contaminated by outliers. Partial quantile regression has also been proposed in \cite{yu2016partial}. However, no theoretical guarantee was originally provided due to the iterative nature of the algorithm and the non-smoothness of quantile loss function. 

 This paper's most important contribution is its proposal of APQR, an alternative formulation of the original PQR method, that guarantees convergence and establishes various asymptotic properties such as the consistency and asymptotic normality of APQR estimators. Analogous  properties for the original PQR method, on the other hand, remain very difficult to study. Analyses of simulated and real data show that our proposed APQR method performs at least comparably to PQR and PLS methods and generally outperforms QRfPC. In particular, APQR can achieve prediction errors competitive with that of other methods, but using a smaller functional basis. This property of APQR, in addition to its asymptotic properties, makes APQR a powerful predictive tool in both theory and practice.

\backmatter

\bmhead{Acknowledgments}

Drs. Linglong Kong, Ivan Mizera, and Bei Jiang acknowledge funding support for this research from the Natural Sciences and Engineering Research Council of Canada (NSERC). Dr. Dengdeng Yu acknowledges funding support for this research from the start up grant from university of Texas at Arlington. Drs. Linglong Kong and Dengdeng Yu further acknowledge funding support from the Canadian Statistical Sciences Institute (CANSSI).

\begin{appendices}

%\section{Section title of first appendix}\label{secA1}

\section{PROOFS AND VERIFICATION}
The appendix contains detailed proofs of the results that are missing in the main paper.
\subsection{Verification of Assumptions of Theorem 1}

%--------------------- Begin of Move to appendix --------------------------------% 
Let us first check assumption (i). For a given fixed $N$, we have that $\rho_{\tau \nu_N}$ is differentiable, and so $l_N(\boldsymbol\upsilon)$ is continuous. As $\lvert\lvert\boldsymbol\upsilon\rvert\rvert \to \infty$, that is, as $\lvert\lvert(\alpha,\boldsymbol\beta)\rvert\rvert \to \infty$ or $\lvert\lvert\mathbf{C}\rvert\rvert \to \infty$, the function $l_N$ should approach $-\infty$. Therefore, $l_N$ is coercive.
To verify assumption (ii), observe that the function $y_i - \alpha - \mathbf{x}_i^T \boldsymbol\beta-
\mathbf{z}_{i}^T \mathbf{C}_{} \mathbf{1}_K $ is an affine function about $\mathbf{C}$. Since $- \rho_\tau$ and its approximation $-\rho_{\tau \nu_N}$ are strictly concave,
we have that $l_N(\boldsymbol\upsilon)$ and $l(\boldsymbol\upsilon)$ are both strictly concave about $\mathbf{C}$. Since they are also strictly concave about $\alpha$, we have strict concavity of $l_N(\boldsymbol\upsilon)$ and $l(\boldsymbol\upsilon)$ about $\boldsymbol\upsilon$. Lastly, assumption (iii) assures that a locally optimized point is isolated. One important property regarding isolated stationary points is that, if the Hessian matrix at a stationary point is nonsingular, then the stationary point is an isolated one (\cite{golubitsky2012stable}). Alternatively, we can impose Condition A1 of \cite{koenker2005quantile} requiring distribution functions $F_i$ to be absolutely continuous with continuous density $f_i$ uniformly bounded away from $0$ and $\infty$ at the point $\xi_i(\tau) = F_i^{-1}(\tau)$, where $F_i$ is the conditional distribution of $y_i$ given $\mathbf{z}_{i}$. Lemma 2 of \cite{hjort2011asymptotics} states that, if we have a sequence of convex functions $l_N(\boldsymbol\upsilon)$ defined on an open convex set $S$ in $\mathbb{R}^{Kd+p+1}$ on which $l_N$ converges pointwisely to $l$, then $\sup_{\boldsymbol\upsilon \in K}\lvert l_N(\boldsymbol\upsilon )-l(\boldsymbol\upsilon )\rvert$ approaches zero for each compact subset $\mathbf K$ of $S$.
As long as the maxima is a unique interior point $S$, we have that the maxima of $l_N$ will approach the maxima of $l$.
%--------------------- End of Move to appendix --------------------------------% 

\subsection{Proof of Theorem 3}

This result can be verified by using Theorem 1 of \cite{rothenberg1971identification}, reproduced below in Lemma \ref{chap2lemrot1971}. 
The regularity assumptions of this lemma are satisfied by the current model since
(1) the parameter space $S$ is open,
(2) the density $p(y,\mathbf{z}\mid\mathbf{C})$ is proper for all $\mathbf{C} \in S$, 
(3) the support of the density $p(y,\mathbf{z}\vert\mathbf{C})$  is the same for all $\mathbf{C} \in S$,
(4) the log density $l_N(\mathbf{C}\vert y,\mathbf{z}) = \ln p(y,\mathbf{z}\vert \mathbf{C})$ is continuously differentiable, and
(5) the information matrix $I_N(\mathbf{C})$ is continuous in $\mathbf{C}$ by \emph{Theorem 2}. Then by Lemma \ref{chap2lemrot1971}, $\mathbf{C}$ is locally identifiable if and only if $I_N(\mathbf{C})$ is nonsingular.

\begin{lemma}[\cite{rothenberg1971identification}, Theorem 1]
\label{chap2lemrot1971}
 Let $\theta_0$ be a regular point of the information matrix $I(\theta)$. Then $\theta_0$ is locally identifiable if and only if $I(\theta_0)$ is nonsingular.
\end{lemma}

\subsection{Proof of Theorem 4}
\begin{proof}
 For simplicity, we omit
$\alpha$ and $\boldsymbol\beta$, though the conclusions
generalize easily to the case with them.
we want to show that the consistency of the estimated factor matrix $\hat{\mathbf{C}}_n$. The following well-known theorem is our major tool for establishing consistency.\\

\begin{lemma}[\cite{van2000asymptotic}, Theorem 5.7]
\label{chap2lem5}%\emph{Lemma 5. (van der Vaart, 1998, Theorem 5.7)}
Let $M_n$ be random functions and let $M$ be a fixed function of $\theta$ such that for every $\epsilon > 0 $
\begin{eqnarray*}
\sup\limits_{\theta \in \nu} \left\vert M_n(\theta) - M(\theta)\right\vert &\to& 0 \quad \textrm{ in probability},\\
\sup\limits_{\theta : d(\theta,\theta_0) \geq \epsilon} M(\theta) &<& M(\theta_0).
\end{eqnarray*}
Then any sequence of estimators $\hat{\theta}_n$ with $M_n(\hat{\theta}) \geq M_n(\theta_0) - o_P(1)$ converges in probability to $\theta_0$.
\end{lemma}

To apply \emph{Lemma \ref{chap2lem5}} in our setting, we take the nonrandom function $M$ to be $\mathbf{C} \mapsto \mathbb{P}_{\mathbf{C}_0} \left[ l_N (Y,Z\mid\mathbf{C})\right]$
and the sequence of random functions to be $M_n: \mathbf{C} \mapsto \frac{1}{n} \sum_{i=1}^n l_N (y_i,z_i\mid\mathbf{C}) = \mathbb{P}_n M$,
where $\mathbb{P}_n$ denotes the empirical measure under $\mathbf{C}_0$.
Then $M_n$ converges to $M$ a.s. by strong law of large number.
The second condition requires that $\mathbf{C}_0$ is a well-separated maximum of $M$.
This is guaranteed by the (global) identifiability of $\mathbf{C}_0$ and information inequality. %{\color{red} information inequality}.
The first uniform convergence condition is most convenient and is verified by the Glivenko-Cantelli theory \citealp{van2000asymptotic}.

The density is $p_{\mathbf{C}} (y \mid \mathbf{z}) = const \cdot \exp\left[ - \rho_{\tau \nu_N} \left(y-\eta(\mathbf{C},\mathbf{z})\right)\right]$
where $\eta(\mathbf{C},\mathbf{z}) = \langle \mathbf{C}, \mathbf{z}\rangle$.
Take $m_{\mathbf{C}}=ln\left[(p_{\mathbf{C}}+p_{\mathbf{C}_0})/2\right]$.
First we show that $\mathbf{C}_0$ is a well-separated maximum of the function $M(\mathbf{C}):=\mathbb{P}_{\mathbf{C}_0} m_{\mathbf{C}}$.
The global identifiability of $\mathbf{C}_0$ and information inequality guarantee that $\mathbf{C}_0$ is the unique maximum of $M$.
To show that it is a well-separated maximum, we need to verify that $M(\mathbf{C}_k) \to M(\mathbf{C}_0)$ implies $\mathbf{C}_k \to \mathbf{C}_0$.

Suppose $M(\mathbf{C}_k) \to M(\mathbf{C}_0)$, then $\langle \mathbf{C}_k, \mathbf{Z} \rangle \to \langle \mathbf{C}_0, \mathbf{Z} \rangle$ in probability. If $\mathbf{C}_k$ are bounded, then $\mathbf{E} \left[ \langle \mathbf{C}_k-\mathbf{C}_0, \mathbf{Z} \rangle^2\right] \to 0$
and $\mathbf{C}_k \to \mathbf{C}_0$ by nonsigularity of $\mathbf{E}\left[ (\mathrm{vec}\mathbf{Z})(\mathrm{vec}\mathbf{Z})^T\right]$.
On the other hand, $\mathbf{C}_k$ can not run to infinity.
If they do, then $\langle \mathbf{C}_k, \mathbf{Z} \rangle / \left\vert \mathbf{C}_k \right\vert  \to 0$ in probability
which in turn implies that $\mathbf{C}_k / \left\vert  C_k\right\vert  \to \mathbf{0}.$ %({\color{red}?})

For the uniform convergence, we see that the class of functions  $\{ \langle \mathbf{C},\mathbf{Z}\rangle, \mathbf{C} \in S \}$
forms a VC class.
This is true because it is a collection of %{\color{red}finite}
number of polynomials of degree 1 and then apply the VC vector space argument
(\cite{van2000weak}, 2.6.15).
This implies that $\{ \eta(\langle \mathbf{C},\mathbf{Z}\rangle), \mathbf{C} \in S \}$ is a VC class since $\eta$ is a monotone function
(\cite{van2000weak}, 2.6.18).

Now $m_{\mathbf{C}}$ is Lipschitz in $\eta$ since
\begin{eqnarray*}
\frac{\partial m_{\mathbf{C}}}{\partial \eta}
= \frac{const \cdot \exp\left[ - \rho_{\tau \nu_N} \left(y-\eta\right)\right] \cdot \rho^\prime_{\tau \nu_N} \left(y-\eta\right)}
{const \cdot \exp\left[ - \rho_{\tau \nu_N} \left(y-\eta\right)\right]+
const \cdot \exp\left[ - \rho_{\tau \nu_N} \left(y-\eta_0\right)\right]}\\
= \frac{ \rho^\prime_{\tau \nu_N} \left(y-\eta\right)}
{1+
\exp\left[ \rho_{\tau \nu_N} \left(y-\eta\right) - \rho_{\tau \nu_N} \left(y-\eta_0\right)\right] }
\leq \sup \left\vert \rho^\prime_{\tau \nu_N} (\cdot) \right\vert =const
\end{eqnarray*}
The last equality holds since $\rho_{\tau \nu_N}(u) \to \rho_\tau(u)$ as $N \to \infty$,
which also implies that $\rho^\prime_{\tau \nu_N}(u) \to \rho^\prime_{\tau}(u)$ as $N \to \infty$ except for $u=0$.
And we know that $\rho^\prime_{\tau \nu_N}(0)=0$.
Similarly we can show that $m_{\mathbf{C}}$ is Lipschitz in $\eta_0$.
A Lipschitz composition of a Donsker class is still a a Donsker class (\cite{van2000asymptotic}, 19.20).
Therefore $\left\{ \mathbf{C} \mapsto m_\mathbf{C} \right\}$
is a bounded Donsker class with the trivial envelope function $1$.
A Donsker class is certainly a Glivenko-Cantelli class.
Finally the Glivenko-Cantelli theorem establishes the uniform convergence condition required by \emph{Lemma \ref{chap2lem5}}.

When the parameter is restricted to a compact set, $\eta(\langle \mathbf{C}, \mathbf{Z}\rangle)$ is confined in a bounded interval
and the $l_N$ is Lipschitz on the finite interval.
It follows that $\{ l_N(\mathbf{C}) = l_N \circ \eta \circ \langle \mathbf{C}, \mathbf{Z}\rangle, \mathbf{C} \in S\}$ is a Donsker class
as composition with a monotone or Lipschitz function preserves the Donsker class.
Therefore the Glivenko-Cantelli theorem establishes the uniform convergence.
Compactness of parameter space implies that $\mathbf{C}_0$ is a well separated maximum if
it is the unique maximizer of $M(\mathbf{C}) = \mathbb{P}_{\mathbf{C}_0} m_{\mathbf{C}}$ (\cite{van2000asymptotic}, Exercise 5.27).
Uniqueness is guaranteed by the information equality whenever $\mathbf{C}_0$ is identifiable.
This verifies the consistency for quantile regression.\\
\end{proof}

\begin{lemma}
\label{chap2lem3}
%\emph{Lemma 3}.
Tensor quantile linear regression model (3) is quadratic mean differentiable (q.m.d.).\\
\end{lemma}
\begin{proof}
By a well-known result (\cite{van2000asymptotic}, Lemma 7.6), it suffices to verify that the density is continuously differentiable
in parameter for $\mu$-almost all $x$ and that the Fisher information matrix exists and it continuous. The derivative of density is
$$\nabla l_N(\mathbf{C}) =  - \sum_{i=1}^{n} \rho^\prime_{\tau \nu_N}\left( \eta_i(\mathbf{C}) \right) \cdot \nabla \eta_i(\mathbf{C}),$$
which is well-defined and continuous by Proposition 2.
The same proposition shows that  the information matrix exists and is continuous.
Therefore the tensor quantile linear regression model is q.m.d.
\end{proof}

\subsection{Proof of Theorem 5}

\begin{proof}
The following result relates asymptotic normality to the density that satisfies q.m.d.\\
\begin{lemma}
\label{chap2lem6}
%{Lemma 6}. \citealp{van2000asymptotic}, Theorem 5.39) suppose that the model $(P_\theta: \theta \in \nu)$ is q.m.d.
At an inner point $\theta_0$ of $\nu \subset \mathbf{R}^k$. Furthermore, suppose that there exists a measurable function
$\dot{l}$ with $\mathbf{P}_{\theta_0} \dot{l}^2 < \infty$ such that, for every $\theta_1$ and $\theta_2$ in a neighbourhood of $\theta_0$,
$$\left\vert  ln p_{\theta_1}(x) - ln p_{\theta_2} (x)\right\vert  \leq \dot{l} (x) \left\vert  \theta_1 -\theta_2 \right\vert .$$
If the Fisher information matrix $I_{\theta_0}$ is nonsingular and $\hat{\theta_n}$ is consistent, then
$$\sqrt{n} (\hat{\theta_n}-\theta_0) = I^{-1}_{\theta_0} \sum_{i=1}^{n} \dot{l}_{\theta_0} (X_i) + o_{P_{\theta_0}}(1).$$
In particular, the sequence $\sqrt{n} (\hat{\theta}_n -\theta_0)$ is asymptotically normal with mean zero and covariance matrix $I^{-1}_{\theta_0}$.
\end{lemma}

\emph{Lemma \ref{chap2lem3}} shows that tensor quantile regression linear model is q.m.d. By \emph{Theorem 2} and chain rule,
the score function
$$\dot{l}_N(\mathbf{C}) = - \sum_{i=1}^{n} \rho^\prime_{\tau \nu_N}\left( \eta_i(\mathbf{C}) \right) \cdot \nabla \eta_i(\mathbf{C})$$
is uniformly bounded in $y$ an $\mathbf{x}$ and continuous in $\mathbf{C}$ for every $y$ and $\mathbf{x}$ with $\mathbf{C}$ ranging over a compact set of $S_0$.
For sufficiently small  neighbourhood U of $S_0$, $\sup_U \left\vert  \dot{l}_{N}(\mathbf{C}) \right\vert $ is square-integrable. Thus the local Lipschitz condition is satisfied and \emph{Lemma \ref{chap2lem6}} applies.
\end{proof}

\end{appendices}

\bibliography{sample}

\end{document}